\documentclass[a4paper, 12pt]{article}
\usepackage{amsmath}
\usepackage{amsfonts}
\usepackage{amssymb} 
\usepackage{amscd}
\usepackage{theorem} 
\usepackage{eufrak} 
\usepackage{euscript} 
\usepackage[latin1]{inputenc}
\usepackage[OT1]{fontenc}

\newcommand{\N} {{\mathbb N}}
\newcommand{\Z} {{\mathbb Z}}

\newcommand{\C} {{\mathbb C}}

\renewcommand{\epsilon} {\varepsilon}
\renewcommand{\phi} {\varphi}

\newcommand{\Br}{\displaybreak[0]}

\newcommand{\sym} {\text{\rm Sym}}

\newcommand{\Id} {\text{\rm Id}}
\newcommand{\ie} {{\em i.e.}\ }
\newcommand{\tensor} {\otimes} 
\newcommand{\bigtensor} {\bigotimes} 

\newcommand{\proof} {{\em Proof.}\ }
\newcommand{\remark}[1] {{\em Remark #1.}\ }
\renewcommand{\leq}{\leqslant}
\renewcommand{\geq}{\geqslant}

\renewcommand{\b} {{\EuFrak b}}
\newcommand{\g} {\EuFrak g}
\newcommand{\n} {{\EuFrak n}}
\newcommand{\Sfrak} {\EuFrak S}
\newcommand{\Bfrak} {\EuFrak B}
\newcommand{\hsl}[1] {{\hat{\EuFrak sl}}(#1)}
\newcommand{\hg} {\hat{\EuFrak g}}
\newcommand{\h} {\hat{\EuFrak h}}
\newcommand{\hh} {{\EuFrak h}}
\newcommand{\D} {\EuScript D}
\renewcommand{\k} {k^-}
\newcommand{\X} {x^+}
\renewcommand{\S} {\text{\rm S}}
\newcommand{\annu}[1] {\text{\rm Annu}_{#1}}
\newcommand{\sh} {\star}
\newcommand{\ctensor} {\hat \otimes} 
 
\newcommand{\cotensor} {\square}
\newcommand{\hpairing}[1]{{\langle #1 \rangle}}

\newcommand{\QQ}{\frac1{q+q^{-1}}}
\newcommand{\QQij}{\frac{\delta_{i,j}}{q+q^{-1}}}
\newcommand{\complete}[1]{{#1}^{\scriptscriptstyle T}}
\newcommand{\qbinom}[3]{{\genfrac{[}{]}{0pt}{1}{#1}{#2}}_{#3}}

\newtheorem{theorem} {Theorem}
\newtheorem{proposition} {Proposition}
\newtheorem{definition} {Definition}
\newtheorem{lemma} {Lemma}

\begin{document}

\title{On quantum shuffle and \\ quantum affine algebras} 
\date{July 2001}
\author{Pascal Grossé}


\maketitle

\begin{center}
\small 
Institut de Recherche Mathématique Avancée,\\
Université Louis Pasteur,\\
7, rue René Descartes, \\
F-67084 Strasbourg Cedex, \\
France\\
email: grosse@math.u-strasbg.fr
\end{center}

\begin{abstract}
  A construction of the quantum affine algebra $U_q(\hg)$ is given in
  two steps. We explain how to obtain the algebra from its positive
  Borel subalgebra $U_q(\b^+)$, using a construction similar to
  Drinfeld's quantum double. Then we show how the positive Borel
  subalgebra can be constructed with quantum shuffles.
\end{abstract}


\section*{Introduction}

Let $\hg$ be an affine Lie algebra over a simple finite dimensional
Lie algebra $\g$.  The quantum affine algebra $U_q(\hg)$ is defined by
generators and relations in the standard Drinfeld-Jimbo presentation.
However, it is well known that the Kac-Moody affine algebra $\hg$ has
a natural realization as a central extension of the loop algebra
$\g\tensor\C[t,t^{-1}]$. Among few attempts to generalize this to the
quantum case, another realization of the quantum affine algebra is
given by Drinfeld in \cite{Drinfeld:NewRealization}. The algebra
structure of $U_q(\hg)$ is given in term of generating series, where
$\hg$ is an untwisted affine Kac-Moody algebra.  However the
Drinfeld-Jimbo coalgebra structure leads to very complicated formulas,
which cannot be expressed in closed form using generating functions. A
new coalgebra structure was given by Drinfeld (in an unpublished
note), with a quite simple formulation.  Ding, Iohara and Frenkel used
this new Drinfeld comultiplication in
\cite{DingFrenkel:IsomorphismRealization} and
\cite{DingIohara:GeneralizationAffineAlgebra}. In the first part of
this work, we show that the algebra structure of $U_q(\hg)$ can be
derived from a construction adapted from Drinfeld's quantum double,
using the comultiplication expressed in closed form in the new
realization.  On the other hand, in his paper
\cite{Nichols:BialgebrasTypeOne}, Nichols examined the structure of
Hopf bimodules. A particular subalgebra of a cotensor construction was
used by Rosso in \cite{Rosso:QuantumShuffles}, where the algebra
structure is completely described with the help of a braiding defining
the quantum shuffle product. By choosing suitable Hopf bimodules,
Rosso showed that the positive Borel subalgebra of the quantized Hopf
algebras $U_q(\g)$ can be obtained using this shuffle construction. In
the second part of this paper, we show that the positive affine
subalgebra $U_q(\b^+)$ can be constructed using a similar method.

{\bf Acknowledgments.} I would like to thank my advisor Marc Rosso for
all the constant encouragements, for having listened to my endless
questions, and for his enlightening responses and suggestions.


\section{A Hopf algebra structure on $U_q(\hg)$}


\subsection{Drinfeld' new realization of $U_q(\hg)$}
\label{sec:gener-drinf-real}

Let $A=(a_{ij})$ be a symmetrizable Cartan matrix corresponding to a
simple Lie algebra $\g$.  Let $\hg$ be the corresponding non twisted
affine Kac-Moody algebra. The normalized symmetric invariant form on
$\h^*$ will be noted $(\cdot|\cdot)$. The set of simple roots in
${\EuFrak h}^*$ is ${\alpha_1,\ldots,\alpha_n}$. If $q$ is a non zero
generic complex number, \ie $q$ is not a root of $1$, we put $q_i =
q^{(\alpha_i|\alpha_i)/2}$ and $q_{ij} = q_i^{a_{ij}} =
q^{(\alpha_i|\alpha_j)}$. Remark that $q_{ij} = q_{ji}$.

Now, let $f_{ij}(t) = (q_{ij}t-1) / (t-q_{ij})$ a complex valued
function. Let $g \in \C[[t]]$ be the formal series $\sum_{n\geq
  0}g^{(ij)}_nt^n$ where the coefficients $g^{(ij)}_n$ are defined
by the Taylor series of $f$ at zero, \ie $f(t) = \sum_{n\geq
  0}g^{(ij)}_n t^n$ for $|t| \ll 1$.

Remark that $f_{ij}(t)f_{ji}(t^{-1}) = 1$. But even by embedding
$\C[[t]]$ in $\C[[t,t^{-1}]]$ in the canonical way, the same relation
does not hold for the formal series $g$, because we cannot have both
$|t| \ll 1$ and $|t| \gg 1$.  By doing this work, we had in mind a
``functional'' point of view, as in
\cite{DingIohara:GeneralizationAffineAlgebra}, where the authors speak
of functional algebras, or in
\cite{LiguoriMintchev:FockRepresentation}, where the authors have
relations such as $x_i(z)x_j(w)=R_{ij}(w,z)x_j(w)x_i(z)$. For such
relations to hold under the permutation of $z$ and $w$, the {\em
  exchange factors} $R_{ij}$ must satisfy the relation
$R_{ij}(z,w)R_{ji}(w,z) = 1$. But as we shall soon see, commutation
relations cannot be stated in such a clean way in our case. However,
we have the two following useful relations in $\C[[t,t^ {-1}]]$:
\begin{align*}
  ( q_{ij} - t)g(t) & = (q_{ij}t-1) \,,\Br\\
  (q_{ij}t - 1)g(t^{-1}) & = (q_{ij} - t) .\Br
\end{align*}

\begin{definition}[\cite{DingIohara:GeneralizationAffineAlgebra}] \label{def:Uq(g)}
  $U_q(\hg)$ is an associative algebra with unit $1$ and generators
  $\{ x^+_{i,n}, x^-_{i,n}, \phi_{i,-k}, \psi_{i,k}, q^{\pm\frac{c}2}
  \mid i=1\ldots n-1, n\in\Z, k\in\N \}$,  satisfying the following
  relations expressed in terms of generating functions in formal
  variables $z$ or $w$:
  \begin{gather}
    q^{\pm\frac{c}2} \, \text{are central and mutually inverse} \,,\label{eq:qiscentral} \Br\\
    \phi_{i,0}{\psi_{i,0}} = {\psi_{i,0}}\phi_{i,0} = 1 \,,\notag \Br\\
    \phi_i(z)\phi_j(w) = \phi_j(w)\phi_i(z) \,,\label{eq:phicommutation} \Br\\
    \psi_i(z)\psi_j(w) = \psi_j(w)\psi_i(z) \,,\notag \Br\\
    g_{ij}(zw^{-1}q^c)\phi_i(z)\psi_j(w) = g_{ij}(zw^{-1}q^{-c})\psi_j(w)\phi_i(z) \,,\Br\\
    \phi_i(z)x^\pm_j(w) = g_{ij}(zw^{-1}q^{\mp\frac{c}2})^{\pm 1}x^\pm_j(w)\phi_i(z) \,,\label{eq:action_phi}\Br\\ 
    \psi_i(z)x^\pm_j(w) = g_{ij}(wz^{-1}q^{\mp\frac{c}2})^{\mp 1}x^\pm_j(w)\psi_i(z) \,,\notag \Br\\ 
    [x^+_i(z),x^-_j(w)] = \frac{\delta_{i,j}}{q-q^{-1}}\left(\delta(zw^{-1}q^{-c})\psi_i(wq^{\frac12 c}) - \delta(zw^{-1}q^c)\phi_i(zq^{\frac12 c})\right) \,,\notag \Br\\
    (z-q_{ij}w)x^\pm_i(z)x^\pm_j(w) =(q_{ij}z - w) x^\pm_j(w)x^\pm_i(z) \,,\label{eq:commutation_x_i} \Br\\
    \sum_{r=0}^{1-a_{ij}}(-1)^r\qbinom{1-a_{ij}}{r}{q_i}\sym_z\left(x^\pm_i(z_1)\cdots x^\pm_i(z_r)x^\pm_j(w)x^\pm_i(z_{r+1})\cdots x^\pm_i(z_{1-a_{ij}})\right) \,,\label{eq:serres}\Br
  \end{gather}
  where $\phi_i(z) = \sum_{k\leq 0}\phi_{i,k}z^{-k}$, $\psi_i(z) =
  \sum_{l\geq 0}\psi_{i,l}z^{-l}$ and $x^\pm_i(z) = \sum_{k\in\Z}
  x^\pm_{i,k}z^{-k}$. The operator $\sym_z$ denotes symmetrization with
  respect to $z_1,\ldots,z_{1-a_{ij}}$, and $\delta$ is the formal distribution
  with support at 1, that is : $\delta(z) = \sum_{k\in\Z}z^k$. As
  usual, $\qbinom{n}{p}{q}$ is the quantum binomial coefficient, with
  $\qbinom{n}{p}{q} = \frac{[n]_q!}{[p]_q![n-p]_q!}$, $[n]_q! =
  [1]_q[2]_q\cdots[n]_q$ and $[n]_q = \frac{q^n-q^{-n}}{q-q^{-1}}$.
\end{definition}

Relations (\ref{eq:serres}) are the so-called Serre relations.

Remark that products such as $\delta(zw^{-1}q^{-c})\psi_i(wq^{\frac12
  c})$ are well defined and can be easily computed.


\subsection{A Hopf algebra structure for $U_q(\hg)$}
 
As the formulas for the coproduct, counit and antipode for this
algebra will involve infinite expressions, we shall make some
topological completion on the underlying vector space and on the
tensor product.  However, exhibiting an inverse system on $U_q(\hg)$
is not straightforward, due to some non-homogeneous relations ({\em
  e.g.} the commutation relation (\ref{eq:action_phi})). Therefore we
shall follow an indirect path.

Let $\overline U_q(\hg)$ be the free algebra with same generators as
$U_q(\hg)$. Now give $\phi_{i,-k}$ and $\psi_{i,k}$ degree $k$ for $k
\geq 0$.  All other elements get degree $0$.  We then extend the
degree on all the elements of the algebra by summation on the
monomials. For $i\geq 0$, let $\overline U_i$ be the ideal of
$\overline U_q(\hg)$ generated by elements of degree greater than $i$.
We then get an inverse system $(\overline U_q(\hg)/\overline
U_i,p_i)$, where $p_i$ is the natural projection $\overline
U_q(\hg)/\overline U_i \rightarrow \overline U_q(\hg)/\overline
U_{i-1}$ obtained using the following diagram, where the rows are
exact sequences:
\begin{equation*}
  \begin{CD}
    0 @>>> \overline U_{i+1} @>>> \overline U_q(\hg) @>>> \overline U_q(\hg)/\overline U_{i+1} @>>> 0 \\
    @. @VVV @| @VV{p_{i+1}}V @. \\
    0 @>>> \overline U_i @>>> \overline U_q(\hg) @>>> \overline U_q(\hg)/\overline U_i @>>> 0 
  \end{CD}
\end{equation*}
The completion of $\overline U_q(\hg)$ is then $\complete{\overline
  U_q(\hg)} = \varprojlim \overline U_q(\hg)/\overline U_i$ ($\scriptstyle T$ stands for
topological), which leads us to the completion of $U_q(\hg)$: let
$\complete{U_q(\hg)} = \complete{\overline U_q(\hg)} / I$, where $I$
is the two-sided ideal generated by the relations in definition
\ref{def:Uq(g)}. There is a canonical injection from $U_q(\hg)$ into
$\complete{U_q(\hg)}$.

In order to complete the tensor product we need a weaker topology (see
remarks below) on $U_q(\g) \tensor U_q(\g)$. As before, we first
consider $\overline U_q(\hg)$.  Then we give $q^{\pm\frac{c}2}$ degree
0. The generators $x^+_{i, \pm k}$, $x^-_{i,\pm k}$, $\phi_{i,-k}$ and
$\psi_{i, k}$ get degree $k$ for $k \geq 0$.  After having extended
the degree the usual way on monomials, we denote by $\overline V_i$
the two sided ideal of $\overline U_q(\hg)$ of elements of degree at
least $i$. Remark that $\overline U_i$ is a strict subset of
$\overline V_i$.  The tensor product $\overline U_q(\hg) \tensor
\overline U_q(\hg)$ can now be completed the following way: we get an
inverse system by setting $\overline S_i = \overline V_i\tensor
\overline V_i$. The inverse limit $\overline U_q(\hg) \ctensor
\overline U_q(\hg) = \varprojlim{}\overline U_q(\hg) \ctensor
\overline U_q(\hg)/\overline S_i$ is then a completion of the usual
tensor product (the symbol $\ctensor$ stands for completed tensor
product).  Now $U_q(\hg) \tensor U_q(\hg)$ can be completed by \[
U_q(\hg) \ctensor U_q(\hg) = \overline U_q(\hg) \ctensor \overline
U_q(\hg) \Big/ \big(\overline U_q(\hg)\tensor I + I \tensor \overline
U_q(\hg)\big).\] Using those completions we can now put a Hopf algebra
structure on $\complete{U_q(\hg)}$, though in a weak sense:

\begin{proposition}
  $\complete{U_q(\hg)}$ is a Hopf algebra for the coproduct $\Delta:
  \complete{U_q(\hg)} \rightarrow U_q(\hg) \ctensor U_q(\hg)$, the
  antipode $\S: \complete{U_q(\hg)} \rightarrow \complete{U_q(\hg)}$ and
  the counit $\epsilon: \complete{U_q(\hg)} \rightarrow \C$ satisfying the
  following relations:
  \begin{align}
    \Delta(q^{\pm\frac{c}2}) & = q^{\pm\frac{c}2} \ctensor q^{\pm\frac{c}2} \,,\notag\Br\\
    \Delta(x^+_i(z)) & = x^+_i(z) \ctensor 1 + \phi_i(zq^{\frac{c_1}2}) \ctensor x^+_i(zq^{c_1}) \,,\label{eq:deltaxplus}\Br\\
    \Delta(x^-_i(z)) & = 1 \ctensor x^-_i(z) + x^-_i(zq^{c_2}) \ctensor \psi_i(zq^{\frac{c_2}2}) \,,\label{eq:deltaxmoins}\Br\\
    \Delta(\phi_i(z)) & = \phi_i(zq^{-\frac{c_2}2}) \ctensor \phi_i(zq^{\frac{c_1}2}) \,,\notag\Br\\
    \Delta(\psi_i(z)) & = \psi_i(zq^{\frac{c_2}2}) \ctensor \psi_i(zq^{-\frac{c_1}2}) \,,\notag\Br\\
    \S(q^{\pm\frac{c}2}) & = q^{\mp\frac{c}2} \,,\notag\Br\\
    \S(x^+_i(z)) & = -\phi_i(zq^{-\frac{c}2})^{-1}x^+_i(zq^{-c}) \,,\label{eq:Sxplus}\Br\\
    \S(x^-_i(z)) & = -x^-_i(zq^{-c})\psi_i(zq^{-\frac{c}2})^{-1} \,,\label{eq:Sxmoins}\Br\\
    \S(\phi_i(z)) & = \phi_i(z)^{-1} \,,\notag\Br\\
    \S(\psi_i(z)) & = \psi_i(z)^{-1} \,,\notag\Br\\
    \epsilon(q^{\pm\frac{c}2}) & = \epsilon(\phi_i(z)) = \epsilon(\psi_i(z)) = 1\,,\notag\Br\\
    \epsilon(x^\pm_i(z)) & = 0 \,,\notag\Br
  \end{align}
  where $c_1 = c\tensor 1$ and $c_2 = 1 \tensor c$.
\end{proposition}

A proof of this proposition can be found in
\cite{DingIohara:GeneralizationAffineAlgebra}. Although the authors of
this paper do not give many details about topological completion, all
the formulas appearing in their proof fit well within our framework.

Remark that the coefficient of $z^n$ in all those formulas involve
finite sums, except in (\ref{eq:deltaxplus}), (\ref{eq:deltaxmoins}),
(\ref{eq:Sxplus}) and (\ref{eq:Sxmoins}). But all infinite sums
converge thanks to the inverse limit topology on $\complete{U_q(\hg)}$
and on the topological tensor product.

We could ask why the topological completions for $U_q(\hg)$ and for the
tensor product are not the same? First, remark that we have to
complete the algebra so that the antipode has a valid definition. The
completion is weak enough to allow us writing expressions such as
$\phi_i(z)x^+_i(z)$ in $\complete{U_q(\hg)}$, but also strong enough to
forbid expressions such as $g_{ij}(zw^{-1})x^+_i(z)x^+_j(w)$. Without
this last obstruction, the Serre relations would have been a trivial
consequence of the commutation relations (\ref{eq:commutation_x_i})
(we shall come back on that later). The topology on the tensor product
is a weaker one, in the sense that expressions such as
$g_{ij}(zw^{-1})x^+_i(z) \ctensor x^+_j(w)$ are allowed, but not
$g_{ij}(zw^{-1})x^+_i(z)x^+_j(w) \ctensor 1$. This discrimination will
be one of the key tool in the sequel. We shall give more details later.
\label{remark-topology}


We would like to obtain the whole algebra $\complete{U_q(\hg)}$ using
Drinfeld's double construction, by giving a Hopf pairing between
suitable ``positive'' and ``negative'' subalgebras of
$\complete{U_q(\hg)}$.  But we shall see in the next part that such a
pairing cannot exist in our case. We will nonetheless exhibit a
construction similar to Drinfeld's one and sharing many properties
(but not all) with it. This can be thought of as a weak quantum double
construction, where we get the algebra structure of $U_q(\hg)$ using a
weak pairing, but of course nothing for the coalgebra structure
because the topology is the main obstruction to the mere existence of
the pairing. Nevertheless, the restriction of the coproduct to the
non-completed subalgebra $U_q(\hg)$ of $\complete{U_q(\hg)}$ coincide
with the formula for the coproduct on the tensor product of the
positive and negative borel subalgebras. Moreover, the pairing will be
a non-degenerate one, which will be of essential importance for the
last part of this work.


\subsection{A weak Hopf pairing between $U_q(\b^+)$ and $U_q(\b^-)$}

Let $\overline U_q(\n^+)$ be the free algebra generated by $x^+_{i,n}$
and $\phi_{i,k}$ for $i=1,\ldots,n-1$, $n\in\Z$ and $k < 0$.  The
upper triangular part of $U_q(\hg)$ is $U_q(\b^+)$ with relations
(\ref{eq:qiscentral}), (\ref{eq:phicommutation}) and
(\ref{eq:action_phi}) given in definition \ref{def:Uq(g)}. The free
algebra with the same basis is $\overline U_q(\b^+) = \overline
U_q(\n^+) \tensor \C[\phi_{i,0}^{\pm 1},q^{\pm\frac12 c}]$.  On the
negative side, we get the algebra $\overline U_q(\n^-)$ generated by
$x^-_{i,m}$ and $\psi_{i,l}$ for $i=1,\ldots,n-1$, $m\in\Z$ and $l >
0$ . The lower triangular part of $U_q(\hg)$ is $U_q(\b^-)$, with the
appropriate relations. The corresponding free algebra is $\overline
U_q(\b^-) = \overline U_q(\n^-) \tensor \C[\psi_{i,0}^{\pm
  1},q^{\pm\frac12 c'}]$, where $c'$ behave as $c$. In the sequel, we
shall also consider the sub algebra $\overline U_q(\hh^+)$ of
$\overline U_q(\b^+)$ generated by $\phi_{i,k}$, $\phi_{i,0}^{-1}$ and
$q^{\pm\frac{c}2}$ for $i=1,\ldots,n-1$ and $k \leq 0$. Its negative
counterpart, $\overline U_q(\hh^-)$ is constructed similarly.  Remark
that $\overline U_q(\hh^+)$ and $\overline U_q(\hh^-)$ are Hopf
algebras, because no topological completion is needed for the
coproduct and the antipode on them. On the contrary, it is necessary to
consider the topological algebras $\complete{\overline U_q(\b^+)}$,
$\complete{\overline U_q(\b^-)}$, $\complete{U_q(\b^+)}$ and
$\complete{U_q(\b^-)}$ if we want to have an Hopf algebra structure.

Let $Q = \bigoplus_{i=1}^{n-1}\Z\alpha_i$ be the root lattice, and
$Q_+ = \bigoplus_{i=1}^{n-1}\N\alpha_i$ be the positive root lattice.
We want to put a $Q$-gradation on the algebras above. If $\alpha=\sum
m_i\alpha_i$ is a root, let us write $\phi_\alpha = \phi_1^{m_1}\cdots
\phi_{n-1}^{m_{n-1}}$ for an expression of type
$\phi_1^{m_1}(z_1)\cdots \phi_{n-1}^{m_{n-1}}(z_N)$ where the formal
variables $z_i$ are all distinct but we do not want to emphasize on
their actual name. 

The gradation is simply obtained by giving $x_i^\pm$ degree
$\pm\alpha_i$, and degree 0 to $\phi_i$ and $\psi_j$, the degree on a
monomial being computed as usual. As a result, we get a $Q$-gradation
on $\overline U_q(\hg)$, $\overline U_q(\b^\pm)$ and $\overline
U_q(\n^\pm)$ (and their corresponding topological algebras). Moreover,
these algebras are direct sums of their subspaces of fixed degree.




Before going further, let us define the main tool of this part:

\begin{definition}[Weak Hopf pairing]
\label{def:weak-hopf-pairing}
Let $A$ and $B$ be algebras over a field $k$, both embedded in
topological Hopf algebras $\complete{A}$ and $\complete{B}$ with invertible
antipodes. A weak Hopf pairing between $A$ and $B$ is a bilinear form
$\hpairing{\cdot,\cdot} : A\times B \rightarrow k$ such that
\begin{enumerate}
\item for all $a\in A$, for all $b,b' \in B$,
  \[ \hpairing{a,bb'} =
  \sum\hpairing{a_{(1)},b}\hpairing{a_{(2)},b'}, \] \label{itemone}
\item for all $a,a'\in A$, for all $b \in B$,
  \[ \hpairing{aa',b} =
  \sum\hpairing{a,b_{(2)}}\hpairing{a',b_{(1)}}, \] \label{itemtwo}
\item for all $a\in A$, for all $b\in B$, 
  \[ \hpairing{a,1_B} = \epsilon(a) \quad \text{and} \quad
  \hpairing{1_A,b} = \epsilon(b), \]
\item for all $a\in A$, for all $b\in B$, 
  \[\hpairing{\S(a),b} = \hpairing{a,\S^{-1}(b)}, \] \label{itemfour}
\end{enumerate}
with the further conditions that in the sums in \ref{itemone}. and
\ref{itemtwo}. (taken over $\Delta_{\complete{A}}(a)$ and
$\Delta_{\complete{B}}(b)$ respectively) there is a finite number of
non zero summands. Similarly, in the expression in \ref{itemfour}.
there is a finite number of non zero terms on both sides of the
equation.
\end{definition}

Remark that the condition of finiteness imposed on the sums and on
the antipode are important, and must be verified for all combinations of
elements of $A$ and $B$. Thus, when we check for a bilinear form to be
a weak Hopf pairing, we must also check that only finite expressions
are involved.

Given such a weak Hopf pairing between two algebras $A$ and $B$, we
can consider their left and right kernel. As in the usual case, we
have the following result:

\begin{proposition}
  Let $I_A$ and $I_B$ be the left and right kernel of the weak
  pairing. Then $I_A$ and $I_B$ are called weak Hopf ideals and satisfy
  \begin{enumerate}
  \item $I_A$ and $I_B$ are two sided ideals,
  \item $I_A$ and $I_B$ are coideals:
    \begin{align*}
      \Delta_{\complete{A}}(I_A) & \subset I_A\ctensor A + A\ctensor I_A, \\
      \Delta_{\complete{B}}(I_B) & \subset I_B\ctensor B + B\ctensor I_B,
    \end{align*}
  \item $I_A$ and $I_B$ are weakly invariant under the antipode:
    \begin{align*}
      \S(I_A) & \subset \complete{I_A} \\
      \S(I_B) & \subset \complete{I_B}, 
    \end{align*}
  \end{enumerate}
where $\complete{I_A}$ (resp. $\complete{I_B}$) is just the usual
topological closure of $I_A$ in $\complete{A}$ (resp. $I_B$ in
$\complete B$).
\end{proposition}

Remark that $\complete{A}/I_A$ is a completion of $A/I_A$, so that we
get naturally an induced weak Hopf pairing between $A/I_A$ and
$B/I_B$. 

Let us consider the topological Hopf algebras $(\complete{\overline
  U_q(\b^+)}, \Delta_+)$ and $(\complete{\overline U_q(\b^-)},
\Delta_-)$. To define a weak Hopf pairing between $\overline U_q(\b^+)$
and $\overline U_q(\b^-)$, we will have to use a slightly generalized
version of the classical result of
\cite{VanDaele:DualPairsHopfAlgebras}. The pairing is defined on the
generators, and by using the fact that the coproduct is a (possibly
infinite) sum of tensor product of generators, we can check that the
pairing can be computed on any two monomials in a finite number of
steps. Incidentally, we shall check that the infinite sums in the
coproduct lead only to finite sums after having applied the pairing on
each of its summands (remember the finiteness conditions in definition
\ref{def:weak-hopf-pairing}).

\begin{proposition}
  There is an unique $Q$-graded weak Hopf pairing between $\overline
  U_q(\b^+) \hookrightarrow (\complete{\overline U_q(\b^+)},\Delta_+)$
  and $\overline U_q(\b^-) \hookrightarrow (\complete{\overline
    U_q(\b^-)},\Delta_-)$ satisfying the following relations :
  \begin{gather*}
    \hpairing{q^{\pm\frac{c}2},q^{\pm\frac{c'}2}} = 1 \,,\\
    \hpairing{q^{\pm\frac{c}2},\psi_{i,k}} = \hpairing{\phi_{i,-k},q^{\pm\frac{c'}2}} = \delta_{k,0} \,,\\
    \hpairing{q^{\pm\frac{c}2},\psi_{i,0}^{-1}} = \hpairing{\phi_{i,0}^{-1},q^{\pm\frac{c'}2}} = 1 \,,\\
    \hpairing{q^{\pm\frac{c}2},x^-_{i,n}} = \hpairing{x^+_{i,n},q^{\pm\frac{c'}2}} = 0 \,,\\
    \hpairing{\phi_{i,l},x^-_{j,n}} = \hpairing{x^+_{j,n},\psi_{i,k}} = 0\,,\\
    \hpairing{\phi_{i,l},\psi_{j,k}} = g^{(ij)}_{k}\delta_{k,-l} \,,\\
    \hpairing{\phi_{i,0}^{-1},\psi_{j,0}} = \hpairing{\phi_{i,0},\psi_{j,0}^{-1}} = (g^{(ij)}_0)^{-1} \,,\\
    \hpairing{x^+_{i,n},x^-_{j,m}} = \frac{-\delta_{i,j}\delta_{n,-m}}{q+q^{-1}} \,,
  \end{gather*}
  for $n,m\in\Z$, $k\geq 0$, $l\leq 0$ and $i=1,\ldots,n-1$.
  
  Moreover, relations (\ref{eq:qiscentral}), (\ref{eq:phicommutation})
  and (\ref{eq:action_phi}) (the action of $\phi_i$ on $x^+_i$)
  are in the left annihilator ideal $\overline I_+$ of this weak
  pairing, with a similar result for the right annihilator ideal
  $\overline I_-$.
\end{proposition}

Looking at all pairing involving $q^{\pm\frac{c}2}$, we see that for
any element $Y \in \overline U_q(\b^-)$ we have
$\hpairing{q^{\pm\frac{c}2}, Y} = \epsilon(Y)$. Thus, for $X$ and $Y$
being any elements of $\overline U_q(\b^+)$ and $\overline U_q(\b^-)$
we have :
\begin{align*}
  \hpairing{Xq^{\pm\frac{c}2},Y} & = \sum\hpairing{X,Y_{(2)}}\hpairing{q^{\pm\frac{c}2}, Y_{(1)}} \\
  & = \sum\hpairing{X,Y_{(2)}}\epsilon(Y_{(1)}) \\
  & = \hpairing{X,Y}.
\end{align*}

Similarly, we can show that $\hpairing{X,q^{\pm\frac{c'}2}Y} =
\hpairing{X,Y}$ for any $X$ and $Y$ in $\overline U_q(\b^+)$ and
$\overline U_q(\b^-)$. Therefore, when involved in pairing
computations, $q^{\pm\frac{c}2}$ is completely irrelevant (\ie can be
considered as the unit 1). The same behavior holds for
$q^{\pm\frac{c'}2}$, so that in all the following computations we shall
omit any reference to those elements in order to simplify the notations.

Before we go to the general existence and uniqueness proof, let us
show that all elements defined by relations (\ref{eq:qiscentral}),
(\ref{eq:phicommutation}), (\ref{eq:action_phi}) leading to
$U_q(\b^+)$ are in the kernel of this weak Hopf pairing (\ie their
pairing with any other element is well defined and null). It is
sufficient to show that each relation $r$ verifies $\hpairing{r,y} =
0$ where $y$ is a generator of $U_q(\b^-)$. We get the general result
by $\hpairing{r, y_1\ldots y_k} =
\sum\prod_{i=1}^k\hpairing{r_{(i)},y_i}$. The sum can of course be
infinite, but each summand is zero because at least one $r_{(i)}$ is
in the Hopf ideal generated by all the relations.

\begin{enumerate}
\item We have already shown that $q^{\pm\frac{c}2}$ and
  $q^{\pm\frac{c'}2}$ are central elements.
\item We have $r_1 = \phi_{i,0}\phi_{i,0}^{-1} - 1 = 0$:
  \begin{enumerate}
  \item$\hpairing{r_1, 1} = 0$, 
  \item$\hpairing{r_1, \psi_{j,k}} = \sum_{r+s=k}\hpairing{\phi_{i,0},\psi_{j,r}}\hpairing{\phi_{i,0}^{-1},\psi_{j,s}} 
      = g^{(ij)}_0(g^{(ij)}_0)^{-1} - 1 = 0$,
  \item$\hpairing{r_1, x^-_{j,m}} = 0$.
  \end{enumerate}
\item We have $r_2 = \phi_{i,n}\phi_{j,m} - \phi_{j,m}\phi_{i,n} = 0$:
  \begin{enumerate}
  \item$\hpairing{r_2, 1} = \delta_{n,0}\delta_{m,0} - \delta_{m,0}\delta_{n,0} = 0$, 
  \item{\setlength{\abovedisplayskip}{-\the\baselineskip}\begin{align*}
       \hpairing{r_2, \psi_{k,o}} & = 
       \smash[t]{\sum_{r+s=o}\Big(\hpairing{\phi_{i,n},\psi_{k,r}}\hpairing{\phi_{j,m},\psi_{k,s}} 
       - \hpairing{\phi_{j,m},\psi_{k,r}}\hpairing{\phi_{i,n},\psi_{k,s}}\Big)} \\
       & =  \sum_{r+s=o}\Big(g^{(ik)}_r\delta_{n,-r}g^{(jk)}_s\delta_{m,-s}       
       - g^{(jk)}_r\delta_{m,-r}g^{(ik)}_s\delta_{n,-s}\Big) \\
       & = \delta_{n+m,-o}\Big(g^{(ik)}_{-n}g^{(jk)}_{-m} - g^{(jk)}_{-m}g^{(ik)}_{-n}\Big) \\
       & = 0,
    \end{align*}}
  \item{\setlength{\abovedisplayskip}{-\the\baselineskip}\begin{align*}
      \hpairing{r_2, x^-_{k,o}} & =
      \hpairing{\phi_{i,n},x^-_{k,o}}\hpairing{\phi_{j,m},1} - \hpairing{\phi_{j,m},x^-_{k,o}}\hpairing{\phi_{i,n},1} \\
      & \quad + \sum_{l\leq 0}\Big(\hpairing{\phi_{i,n},\psi_{k,l}}\hpairing{\phi_{j,m},x^-_{k,o-l}} - \hpairing{\phi_{j,m},\psi_{k,l}}\hpairing{\phi_{i,n},x^-_{k,o-l}} \Big) \\
      & = 0.
    \end{align*}}
  \end{enumerate}
\item We have $\phi_i(z)x^+_j(w) =
  g_{ij}(zw^{-1}q^{-\frac{c}2})x^+_j(w)\phi_i(z)$, \ie in term of
  underlying generators: $r_3 = \phi_{i,n}x^+_{j,m} - \sum_{0\leq o\leq
    -n}g^{(ij)}_o q^{-o\frac{c}2} x^+_{j,m-o}\phi_{i,n+o} = 0$:
  \begin{enumerate}
    \item $\hpairing{r_3,1} = 0$,
    \item{\setlength{\abovedisplayskip}{-\the\baselineskip}\begin{align*}
          \hpairing{r_3,\psi_{k,p}} & = \smash[t]{\sum_{r+s=p}\Bigg(}\hpairing{\phi_{i,n},\psi_{k,r}}\hpairing{x^+_{j,m},\psi_{k,s}} \\
          & \qquad - \sum_{0\leq o \leq -n}g^{(ij)}_o\hpairing{x^+_{j,m-o},\psi_{k,r}}\hpairing{\phi_{i,n+o},\psi_{k,s}}\Bigg) = 0,
        \end{align*}}
    \item{\setlength{\abovedisplayskip}{-\the\baselineskip}\begin{align*}
          \hpairing{r_3,x^-_{k,o}} & = \hpairing{\phi_{i,n},x^-_{k,o}}\hpairing{x^+_{j,m},1} \\
          & \qquad - \sum_{0\leq p\leq -n}g^{(ij)}_p\hpairing{x^+_{j,m-p},x^-_{k,o}}\hpairing{\phi_{i,n+p},1} \\
          & \quad + \sum_{l\geq 0}\Bigg( \hpairing{\phi_{i,n},\psi_{k,l}}\hpairing{x^+_{j,m},x^-_{k,o-l}} \\
          & \qquad - \sum_{0\leq p\leq -n}g^{(ij)}_p\hpairing{x^+_{j,m-p},\psi_{k,l}}\hpairing{\phi_{i,n+p},x^-_{k,o-l}}\Bigg) \\
          & = \QQ g^{(ij)}_p\delta_{k,j}\delta_{m-p,-o}\delta_{n+p,0} - 
          \QQ g^{(ia)}_l\delta_{k,j}\delta_{l,-n}\delta_{m,l-k} \\
          & = \QQ g^{(ij)}_{-n}\delta_{k,j}\delta_{m+n,-o} -
          \QQ g^{(ij)}_{-n}\delta_{k,j}\delta_{m+n,-o} \\
          & = 0.
        \end{align*}}
  \end{enumerate}
\end{enumerate}

Of course the same computations can be made for the corresponding
relations in $\overline U_q(\b^-)$. 

In order to show the existence of the weak pairing, we have to prove
that it is well defined on any two monomials, as was done in
\cite{VanDaele:DualPairsHopfAlgebras}. But we have to be careful with
comultiplications involving infinite sums. 

A first step towards that goal is given by the following result: given
two monomials $X = x^+_{i_1,n_1}\cdots x^+_{i_k,n_k}$ and $Y =
x^-_{j_1,m_1}\cdots x^-_{j_l,m_l}$, $\Phi$ a monomial in $\overline
U_q(\hh^+)$ and $\Psi$ a monomial in $\overline U_q(\hh^-)$, we have
\begin{align*}
  \hpairing{X\Phi, Y\Psi} = \hpairing{X,Y}\hpairing{\Phi,\Psi}
\end{align*}

Indeed, we have $\Delta_+(X) = X \ctensor 1 + R$, where $R$ is an
(infinite) sum of elements having at least one $x^+$ generator in the
right side of the tensor product. For this reason we have
$\hpairing{R\Delta_+(\Phi),Y\tensor \Psi} = 0$. So we get
\[ \hpairing{X\Phi, Y\Psi} = \sum\hpairing{X\Phi_{(1)},Y}
\hpairing{\Phi_{(2)},\Psi}, \] where this sum is taken over
$\Delta_+(\Phi)$ and so is finite. Using the
same argument with $\Delta_-$, we get:
\begin{align*}
  \hpairing{X\Phi, Y\Psi} & = \sum\hpairing{X\Phi_{(1)},Y}\hpairing{\Phi_{(2)},\Psi}\\
  & = \sum\hpairing{X,Y}\hpairing{\Phi_{(1)},1}\hpairing{\Phi_{(2)},\Psi} \\
  & = \sum\hpairing{X,Y}\epsilon(\Phi_{(1)})\hpairing{\Phi_{(2)},\Psi} \\
  & = \hpairing{X,Y}\hpairing{\Phi,\Psi}.
\end{align*}
Now, the second pair involves only finite sums, so it can be computed using
weak Hopf pairing properties in a finite number of operations. For the
first one, we have to be more careful:

Let $\alpha$ be the degree of $X$ and $\beta$ be the degree of $Y$. We
have  
\begin{align*}
  \Delta_+^{(n)}(x^+_{i_l,n_l}) & = x^+_{i_l,n_l} \ctensor 1 \cdots \ctensor 1 +
  \sum_{k_1\leq 0} \phi_{i_l,k_1}q^{?c}\ctensor x^+_{i_l,n_l-k_1} \ctensor 1
  \cdots \ctensor 1 + \ldots + \\ 
  & \sum_{\substack{k_{n-1}\leq 0\\r_{n-1,1}+\ldots+r_{n-1,n-1} =
    k_{n-1}}}\phi_{i_l,r_{n-1,1}}q^{?c}\ctensor \cdots \ctensor
  \phi_{i_l,r_{n-1,n-1}}q^{?c} \ctensor x^+_{i_l,n_l-k_{n-1}}, 
\end{align*}
We do not care for the coefficients in $q^{?c}$, because those are
irrelevant as regards weak Hopf pairings. Then, it appears that in
\[\hpairing{X,Y} = \hpairing{\Delta^{(k-1)}(X),
  x^-_{j_1,m_1}\ctensor\cdots\ctensor x^-_{j_k,m_k}},\] the only
non-zero term can be those where each $x^+_{i,n}$ is paired with a
corresponding $x^-_{i,m}$. Therefore, let $\sigma \in \Sfrak_k$ be a
permutation satisfying $i_{\sigma(l)} = j_l$ for all $1 \leq l \leq
k$. The non-zero terms in the weak Hopf pairing appear when
$x^-_{j_1,m_1}\ctensor\cdots\ctensor x^-_{j_l,m_l}$ is paired with
\[\prod_{l=1}^k \sum_{k_l\leq 0,
    r_{l,1} + \cdots + r_{l, \sigma(l)-1} = 
  k_l} \phi_{i_l,r_{l,1}}q^{?c}
\ctensor\cdots\ctensor \phi_{i_l,r_{l,\sigma(l)-1}}q^{?c}
\ctensor x^+_{i_l,n_l-k_l} \ctensor 1
\ctensor\cdots\ctensor 1.\] 
As we want to rearrange this expression to compute all the pairings by
applying relations given in definition \ref{def:weak-hopf-pairing}, we
set $\mu = \sigma^{-1}$ in $\Sfrak_k$. The last expression then
becomes : 
\begin{gather*}
\begin{split}
  \smash{\sum_{\substack{1\leq l \leq k,\; k_l \leq 0, \\
      r_{l,1} + \cdots + r_{l,\sigma(l)-1} = k_l}} \bigtensor_{p=1}^k
  \Bigg(\prod_{\substack{p < m \leq k \\ \mu(m) < \mu(p)}}}
    \phi_{\mu(m),r_{\mu(m),p}} & \cdot x^+_{i_{\mu(p)},n_{\mu(p)} -
      k_{\mu(p)}} \\ 
  & \cdot \prod_{\substack{p < m \leq k \\ \mu(m) >
        \mu(p)}} \phi_{\mu(m),r_{\mu(m),p}} q^{?c}\Bigg).
\end{split}
\end{gather*}
This sum still involves an infinite number of terms. But remark that
in expressions of the form $\phi_{i_{n_1},?}\cdot\phi_{i_{n_r},?}
x^+_{i_l,?}  \phi_{i_{n_{r+1}},?}\cdot\phi_{i_{n_p},?}$, we have
$\sigma(n_i) > \sigma(l)$ for all $1\leq i \leq p$. Now, it is easy to
see that a pairing of the form $\hpairing{\phi_{i_1,n_1}
  \cdots\phi_{i_r,n_r}x^+_{i,n}\phi_{i_{r+1},n_{r+1}}\cdots
  \phi_{i_p,n_p}, x^-_{j,m}}$ is always zero but for $i=j$,
$n_{r+1}+\cdots+ n_p = 0$ and $n = m - (n_1+\cdots+n_r)$.  Starting
with $l=k$, we get $n_{\mu(k)} - k_{\mu(k)} = m_k$ for the only
non-zero pairing. Thus, $k_{\mu(k)}$ is fixed and the $r_{k,i}$ for $1
\leq i \leq k$ range in a finite number of values. Thereafter for
$l=k-1$, we see that $k_{\mu(k-1)}$ can take only a finite number of
values (if we want the pairing to be non-zero), so that $r_{k-1,i}$
for $1 \leq i \leq k-1$ range over a finite domain also. It is now
easy to continue backward until $l=1$ and conclude that the pairing is
non-zero only for a finite number of terms in the above sum. Thus, the
pairing of any two expression is computable in a finite number of
steps, and the result is well defined, as stated in the proposition.

We have proved that $\hpairing{X\Phi,Y\Psi}$ is well defined and
unique. Because of relation~(\ref{eq:action_phi}), the pairing between
any monomials $\hpairing{X',Y'}$ can be put under the above form in a
finite number of steps, thus finishing the proof.


\subsection{The weak quantum double $\D(U_q(\b^+),U_q(\b^-))$}

Using our weak Hopf pairing, we get:

\begin{proposition}[Weak quantum double]
There is an algebra structure on $\D(\overline U_q(\b^+), \overline
U_q(\b^-)) = \overline U_q(\b^+)\tensor \overline U_q(\b^-)$ where 
\begin{align}
  (a\tensor 1)(1\tensor b) & = a \tensor b\,,\notag\\
  (1\tensor b)(a\tensor 1) & = \sum\hpairing{a_{(1)}, \S(b_{(1)})}\hpairing{a_{(3)}, b_{(3)}}a_{(2)}\tensor b_{(2)}\,,\label{prop:skew-product}\\
  (a\tensor 1)(a'\tensor 1) & = aa' \tensor 1\,,\notag\\
  (1\tensor b)(1\tensor b') & = 1 \tensor bb'\,,\notag
\end{align}
with unit $1\tensor 1$. Moreover, we have natural embeddings
\begin{alignat*}{2}
  \overline U_q(\b^+) & \longrightarrow \D(\overline U_q(\b^+),\overline U_q(\b^-)) & \qquad \overline U_q(\b^-) & \longrightarrow \D(\overline U_q(\b^+), \overline U_q(\b^-)) \\
  a & \longmapsto a\tensor 1 & \qquad b & \longmapsto 1\tensor b
\end{alignat*}
which are algebra morphisms.
\end{proposition}

The proof (mainly the associativity of the multiplication) is similar
to the non-weak case, so it will not be developed here. Though,
because we are using a weak Hopf pairing, we need to be careful about
the sum appearing in the multiplication (\ref{prop:skew-product}). We
will postpone this verification until the proof of next
proposition.

Remark that the coefficient $-1/(q+q^{-1})$ for the Hopf pairing
between $x^+_{i,n}$ and $x^-_{j,n}$ is here so that we get the right
commutation relations in the weak quantum double. In fact, we have :

\begin{proposition}
  In $\D(\overline U_q(\b^+),\overline U_q(\b^-))$ the following
  commutation relations hold:
  \begin{gather}
    q^{\pm\frac{c}2} \, \text{and} \, q^{\pm\frac{c'}2} \, \text{are central} \,,\notag \Br\\
    g_{ij}(zw^{-1}q^{\frac{c}2\tensor\frac{c'}2})\phi_i(z)\psi_j(w) = g_{ij}(zw^{-1}q^{-\frac{c}2\tensor \frac{c'}2})\psi_j(w)\phi_i(z) \,,\label{eq:commutationphipsi}\Br\\
    \phi_i(z)x^-_j(w) = g_{ij}(zw^{-1}q^{\frac{c}2})^{-1}x^-_j(w)\phi_i(z) \,,\notag\Br\\ 
    \psi_i(z)x^+_j(w) = g_{ij}(wz^{-1}q^{-\frac{c'}2})^{-1}x^+_j(w)\psi_i(z) \,,\notag \Br\\ 
    [x^+_i(z),x^-_j(w)] = \frac{\delta_{i,j}}{q-q^{-1}}\left(\delta(zw^{-1}q^{-c'})\psi_i(wq^{\frac12 c'}) - \delta(zw^{-1}q^c)\phi_i(zq^{\frac12 c})\right) \,.\notag \Br
  \end{gather}
\end{proposition}

\proof The verifications are quite long, but straightforward. We have
\begin{enumerate}
\item $q^{\frac{c}2}$ is central (same demonstration for
  $q^{\frac{c'}2}$):
  \begin{enumerate}
  \item $(1\tensor q^{\frac{c'}2})(q^{\frac{c}2}\tensor 1) =
    q^{\frac{c}2}\tensor q^{\frac{c'}2}$,
  \item{\setlength{\abovedisplayskip}{-\the\baselineskip}\begin{align*}
        \smash[t]{(1\tensor\psi_{j,k})(q^{\frac{c}2}\tensor 1)} & = \smash[t]{\sum_{u+v+w=k}
        \hpairing{1,\S(\psi_{j,u})}\hpairing{1,\psi_{j,w}} q^{\frac{c}2}
        \tensor \psi_{j,v}q^{(w-v)\frac{c'}2}} \\
        & = q^{\frac{c}2} \tensor\psi_{j,k},
      \end{align*}}
  \item We have 
    \begin{align*}
      \Delta^2_-(x^-_{j,m}) & = 1\ctensor 1\ctensor x^-_{j,m} +
      \sum_{k\leq 0}1\ctensor x^-_{j,m-k}\ctensor\psi_{j,k}
      q^{-(m-\frac{k}2)c'} \\
      &\quad + \sum_{\substack{k\leq
          0\\r+s=k}}x^-_{j,m-k}\ctensor\psi_{j,r}q^{-(m-s-\frac{r}2}c'\ctensor\psi_{j,s}q^{-(m-\frac{s}2)c'},
    \end{align*}
    so we get
    \begin{align*}
      (1\tensor x^-_{j,m})(q^{\frac{c}2}\tensor 1) & =
      \hpairing{1,1}\hpairing{1,x^-_{j,m}}q^{\frac{c}2}\tensor 1 + \sum_{k\leq
        0}\hpairing{1,1}\hpairing{1,\psi_{j,k}}q^{\frac{c}2}\tensor
      x^-_{j,m-k} \\
      &\quad + \sum_{\substack{k\leq 0\\r+s=k}}
      \hpairing{1,\S(x^-_{j,m-k})}\hpairing{1,\psi_{j,s}}q^{\frac{c}2}\tensor \psi_{j,r}q^{-(m-s-\frac{r}2)c'} \\
      & = 0 + \delta_{k,0}q^{\frac{c}2}\tensor x^-_{j,m-k} + 0 \\
      & = q^{\frac{c}2}\tensor x^-_{j,m}.
    \end{align*}
  \end{enumerate}
\item Relation
  $g_{ij}(zw^{-1}q^{\frac{c}2\tensor\frac{c'}2})\phi_i(z)\psi_j(w) =
  g_{ij}(zw^{-1}q^{-\frac{c}2\tensor
    \frac{c'}2})\psi_j(w)\phi_i(z)$. We have
  \begin{align*}
    \Delta^2_+(\phi_{i,n}) & = \sum_{r+s+t = n}\phi_{i,r}q^{-(s+t)\frac{c}2}\ctensor\phi_{i,s}q^{(r-t)\frac{c}2}\ctensor\phi_{i,t}q^{(r+s)\frac{c}2} \,,\\
    \Delta^2_-(\psi_{j,m}) & = \sum_{u+v+w = m}\psi_{j,u}q^{(v+w)\frac{c'}2}\ctensor\psi_{j,v}q^{(w-u)\frac{c'}2}\ctensor\psi_{j,w}q^{-(u+v)\frac{c'}2} \,,
  \end{align*}
  so that we get
  \begin{equation*}
    \begin{split}
      (1\ctensor\psi_{j,m})(\phi_{i,n}\ctensor 1) = \smash{\sum_{\substack{r+s+t=n\\
            u+v+w=m}}}\hpairing{\phi_{i,r}, \S(\psi_{j,u})}
      & \hpairing{\phi_{i,t},\psi_{j,w}} \\
      & \qquad \phi_{i,s}q^{(r-t)\frac{c}2}\ctensor
      \psi_{j,v}q^{(w-u)\frac{c'}2}, 
    \end{split}
  \end{equation*}
  where we set $\S(\psi_j(z)) = \psi_j(z)^{-1} = \sum_{k\geq
    0}\S(\psi_{j,k})z^{-k}$. Now, it is easy to see (using suitable
  properties of weak Hopf pairings) that 
  \[\hpairing{\phi_{i,t}q^{(r+s)\frac{c}2}, \psi_{j,w}
    q^{-(u+v)\frac{c'}2}} = \hpairing{\phi_{i,t},\psi_{j,w}} =
  \delta_{w,-t}g^{(ij)}_w.\] 
  For the other pair, we proceed as follows: $\S(\psi_{j,u}) =
  P(\psi_{j,0}^{-1},\psi_{j,0},\ldots,\psi_{j,u})$, where $P$ is a
  polynomial such that each of its monomial $\psi_{j,n_1}^{m_1}\cdots
  \psi_{j,n_k}^{m_k}$ satisfies $\sum_p m_pn_p = u$. So we get exactly
  \begin{align*}
    \hpairing{\phi_{i,r}q^{-(s+t)\frac{c}2}, \S(\psi_{j,u})
      q^{-(v+w)\frac{c'}2}} & =
    \hpairing{\phi_{i,r}, P(\psi_{j,0}^{-1},\psi_{j,0},\ldots,
      \psi_{j,u})} \\ 
    & =  P\left( \hpairing{\phi_{i,r}, \psi_{j,0}^{-1}},
      \hpairing{\phi_{i,r},\psi_{j,0}},\ldots, \hpairing{\phi_{i,r},
        \psi_{j,u}}\right) \\  
    & =  P( {g^{(ij)}_0}^{-1}, g^{(ij)}_0, \ldots, g^{(ij)}_u)
    \delta_{r,-u}\\ 
    & = {g^{(ij)}_u}^{-1}\delta_{r,-u}.
  \end{align*}
  Combining those two results, we get
  \begin{align*}
    (1\ctensor\psi_{j,m})(\phi_{i,n}\ctensor 1) & =
    \sum_{\substack{k\geq 0 \\r+s=k}}\Big({g^{(ij)}_r}^{-1}g^{(ij)}_sq^{-(a+c)\frac{c}2}\ctensor
    q^{-(a+c)\frac{c'}2} \Big)\phi_{i,n+k}\ctensor \psi_{j,m-k} \\
    & = \sum_{k\geq 0} h^{(ij)}_k \phi_{i,n+k}\ctensor \psi_{j,m-k},
  \end{align*}
  where $h_{ij}(t) = \sum_{k\geq 0}h^{(ij)}_kt^k = g_{ij}(tq^{\frac{c}2\tensor\frac{c'}2})g_{ij}(tq^{-\frac{c}2\tensor\frac{c'}2})^{-1}$.
  Using currents, the last statement gives
  \[\psi_j(w)\phi_i(z) = h(zw^{-1})\phi_i(z)\psi_j(w),\] which is
  exactly what we wanted to show.
\item Relation $x^-_j(w)\phi_i(z) =
  g_{ij}(zw^{-1}q^{\frac{c}2})\phi_i(z)x^-_j(w)$, \ie \[(1\tensor
  x^-_{j,m})(\phi_{i,n}\tensor 1) = \sum_{0\leq k\leq
    -n}g^{(ij)}_kq^{k\frac{c}2}\phi_{i,n+k}\tensor x^-_{j,m-k}.\]
  We have
  \begin{align*}
    (1\tensor x^-_{j,m})& (\phi_{i,n}\tensor 1) \\
    & = \sum_{r+s+t=n}\hpairing{\phi_{i,r},\S(1)} \hpairing{\phi_{i,t},x^-_{j,m}} \phi_{i,s}q^{(r-t)\frac{c}2}\tensor 1 \\
    &\quad + \sum_{\substack{k\leq 0\\r+s+t=n}}\hpairing{\phi_{i,r},\S(1)}\hpairing{\phi_{i,t},\psi_{j,k}}\phi_{i,s}q^{(r-t)\frac{c}2}\tensor x^-_{j,m-k}\\
    &\quad + \smash[b]{\sum_{\substack{k\leq 0,u+v=k\\r+s+t=n}}}\hpairing{\phi_{i,r},\S(x^-_{j,m-k})}\hpairing{\phi_{i,t},\psi_{j,v}} \\
    & \qquad\qquad\qquad\qquad\qquad \phi_{i,s}q^{(r-t)\frac{c}2}\tensor \psi_{j,u}q^{-(m-v-\frac{u}2)c'}\\
    & = 0 + \sum_{\substack{k\leq 0\\r+s+t=n}}\delta_{r,0}g^{(ij)}_k\delta_{t,-k}\phi_{i,s}q^{(r-t)\frac{c}2}\tensor x^-_{j,m-k} + 0 \\
    & = \sum_{0\leq k\leq -n}g^{(ij)}_k\phi_{i,n+k}q^{k\frac{c}2}\tensor x^-_{j,m-k}.
  \end{align*}
\item Relation{\setlength{\abovedisplayskip}{-\the\baselineskip}\begin{align*}
      [x^+_i(z),x^-_j(w)] & = \smash[t]{\frac{\delta_{i,j}}{q-q^{-1}}\Big(\delta(zw^{-1}q^{-c'})\psi_i(wq^{\frac12 c'})} \\
      & \qquad\qquad - \delta(zw^{-1}q^c)\phi_i(zq^{\frac12 c})\Big)
    \end{align*}}
  \ie 
  \begin{align*}
    (1\tensor x^-_{j,m})(x^+_{i,n}\tensor 1) & = x^+_{i,n}\tensor x^-_{j,m} + \QQij 1\tensor\psi_{j,n+m}q^{\frac{n-m}2c'} \\
    & \qquad\qquad - \QQij \phi_{i,n+m}q^{\frac{m-n}2c}\tensor 1.
  \end{align*}
  We have
  \begin{align*}
    (1\tensor x^-_{j,m})&(x^+_{i,n}\tensor 1) = \hpairing{x^+_{i,n},1}\hpairing{1,x^-_{j,m}}1\tensor 1 \\
    &\quad + \sum_{k\geq 0}\hpairing{x^+_{i,n},1}\hpairing{1,\psi_{j,k}}1\tensor x^-_{j,m-k} \\
    &\quad + \sum_{\substack{k\geq 0\\r+s=k}}\hpairing{x^+_{i,n}, \S(x^-_{j,m-k})}\hpairing{1,\psi_{j,s}}1\tensor \psi_{j,r}q^{-(m-s-\frac{r}2)c'}\\
    &\quad + \sum_{l\leq 0}\hpairing{\phi_{i,l},1}\hpairing{1,x^-_{j,m}}x^+_{i,n-l}\tensor 1 \\
    &\quad + \sum_{\substack{l\leq 0\\k\geq 0}}\hpairing{\phi_{i,l},1}\hpairing{1,\psi_{j,k}}x^+_{i,n-l}\tensor x^-_{j,m-k}\\
    &\quad + \smash[b]{\sum_{\substack{l\leq 0\\k\geq0, r+s=k}}}\hpairing{\phi_{i,l},S(x^-_{j,m-k})}\hpairing{1,\psi_{j,s}}\\
    &\qquad\qquad\qquad\qquad x^+_{i,n-l}\tensor \psi_{j,r}q^{-(m-s-\frac{r}2)c'}\\
    &\quad + \sum_{\substack{l\leq 0\\t+u=l}}\hpairing{\phi_{i,t},1}\hpairing{x^+_{i,n-l},x^-_{j,m}}\phi_{i,u}q^{-(n-t-\frac{u}2)c}\tensor 1\\
    &\quad + \sum_{\substack{l\leq 0,t+u=l\\k\geq 0}}\hpairing{\phi_{i,t},1}\hpairing{x^+_{i,n-l},\psi_{j,k}}\phi_{i,u}q^{-(n-t-\frac{u}2)c}\tensor x^-_{j,m-k}\\
    &\quad + \smash[b]{\sum_{\substack{l\leq 0,t+u=l\\k\geq 0,r+s=k}}}\hpairing{\phi_{i,t},\S(x^-_{j,m-k})}\hpairing{x^+_{i,n-l},\psi_{j,s}}\\
    &\qquad\qquad\qquad\qquad\phi_{i,u}q^{-(n-t-\frac{u}2)c}\tensor \psi_{j,r}q^{-(m-s-\frac{r}2)c'}.
  \end{align*}
  Among the 9 summands of the right hand side of this equation, it
  is trivial to see that the first, second, fourth and eighth are
  0. For the remaining summands, remark that 
  \[\S(x^-_{j,m-k}) = -\sum_{p\geq 0}
  x^-_{j,m-k-p}\S(\psi_{j,p})q^{(m-k-\frac{p}2)c'},\]
  so that $\hpairing{\phi_{i,\cdot},\S(x^-_{j,m-k})} = 0$ and
  $\hpairing{x^+_{i,n},\S(x^-_{j,m-k})} = \QQij \delta_{-n,m-k}$.
  This shows that the sixth and ninth terms are also zero. Now, the
  fifth term is exactly $x^+_{i,n}\tensor x^-_{j,m}$, the third term
  is \[\sum_{\substack{k\geq
      0\\r+s=k}}\QQij\delta_{-n,m-k}\delta_{s,0}1\tensor
  \psi_{j,k}q^{-(m-\frac{k}2)c'},\]
  which is equal to $\QQij1\tensor
  \psi_{j,n+m}q^{\frac{n-m}2c'}$. Finally, the seventh and last term
  is \[-\sum_{\substack{l\leq
      0\\t+u=l}}\QQij\delta_{n-l,-m}\delta_{t,0}
  \phi_{i,l}q^{-(n-\frac{l}2)c}\tensor 1,\]
  which is equal to $-\QQij\phi_{i,n+m}q^{\frac{m-n}2c}\tensor 1$,
  thus achieving the proof.
\end{enumerate}

Remark that this proof also shows the multiplication on the weak
Quantum double is well defined (all sums in the product are finite).

Remember that relations (\ref{eq:qiscentral}), (\ref{eq:phicommutation}), 
(\ref{eq:action_phi}) and their negative counterparts are in the
kernel of this weak Hopf pairing. A natural question arising is whether
the weak Hopf ideal they generate are the whole annihilator
ideals. Actually this is not the case, and we have the following result:

\begin{proposition} \label{prop:RelationsAnnihilator}
  The elements
  \begin{gather}
    (z-q_{ij}w)x^+_i(z)x^+_j(w) - (q_{ij}z-w)x^+_j(w)x^+_i(z) \label{eq:commuteannihilator}
  \end{gather}
  are in the annihilator ideal $\overline I_+$ of the weak Hopf
  pairing between $\overline U_q(\b^+)$ and $\overline U_q(\b^-)$.
  Moreover, the two-sided ideal they generate is a weak Hopf ideal. A
  similar statement holds for the $x^-_i(z)$ ($i=1,\ldots,n-1$) and
  the corresponding annihilator ideal $I_-$.
\end{proposition}

\proof For the sake of this proof, we will denote by $I^0_+$ and
$I^0_-$ the weak Hopf ideals generated by the relations in proposition
(\ref{prop:skew-product}). We can now work in the algebra
$\D(\overline U_q(\b^+)/I^0_+,\overline U_q(\b^-)/I^0_-)$ with the
induced weak Hopf pairing and the induced $Q$-gradation (the relation
defining $I^0_+$ and $I^0_-$ are $Q$-homogeneous).

Recall that a quasi-primitive element of $\overline U_q(\b^+)/I^0_+$
is an element $x$ such that $\Delta_+(x) = x\ctensor h + h'\ctensor
x$, where $h$ and $h'$ are in the subalgebra generated by $\phi_i(z)$
for $i=1,\ldots,n-1$ and $q^{\pm\frac{c}2}$, {\em i.e.}  $U_q(\h^+)$.
It now easy to check that (\ref{eq:commuteannihilator}) is
quasi-primitive : 
\begin{align*}
& \Delta_+\left((z-q_{ij}w)x^+_i(z)x^+_j(w) -
  (q_{ij}z-w)x^+_j(w)x^+_i(z)\right) \\
& \qquad = \left[(z-q_{ij}w)x^+_i(z)x^+_j(w) -
  (q_{ij}z-w)x^+_j(w)x^+_i(z)\right] \ctensor 1 \\
& \qquad \phantom{=} +
  \phi_i(zq^{\frac{c_1}2})\phi_j(wq^{\frac{c_1}2}) \ctensor
  \Big[(z-q_{ij}w)x^+_i(zq^{c_1})x^+_j(wq^{c_1}) \\
& \qquad\phantom{= + \phi_i(zq^{\frac{c_1}2})\phi_j(wq^{\frac{c_1}2}) \ctensor
  \Big[}
  - (q_{ij}z-w)x^+_j(wq^{c_1})x^+_i(zq^{c_1})\Big] \\
& \qquad\phantom{=} + 
  \left[(z-q_{ij}w)\phi_i(zq^{\frac{c_1}2})x^+_j(w) - 
  (q_{ij}z-w)x^+_j(w)\phi_i(zq^{\frac{c_1}2})\right] \ctensor
  x^+_i(zq^{c_1})\\
& \qquad\phantom{=}+
  \left[(z-q_{ij}w)x^+_i(z)\phi_j(wq^{\frac{c_1}2}) - 
  (q_{ij}z-w)\phi_j(wq^{\frac{c_1}2})x^+_i(z)\right] \ctensor
  x^+_j(wq^{c_1})
\end{align*}
The last two elements of this sum are zero thanks to the commutation
relations in $\overline U_q(\b^+)/I^0_+$. Now, quasi-primitive
elements of $\overline U_q(\b^+)/I^0_+$ are orthogonal to decomposable
elements of $\overline U_q(\b^-)/I^0_-$. Remark that all elements of 
degree $-\alpha_i-\alpha_j$ are decomposable. Using the fact that the
weak Hopf pairing is $Q$-graded, we get the first result. 

We can then consider the two sided ideal $I^1_+$ of $\overline
U_q(\b^+)$ generated by $I^0_+$ and the above relation. Likewise, we
have a two sided ideal $I^1_-$ on the negative side. Those two ideals
are again weak Hopf ideals.

Unlike the classical case, the Serre relations are not
quasi-primitive in $\overline U_q(\b^+)$ and $\overline
U_q(\b^-)$. But they are quasi-primitive modulo commutation relations
(\ref{eq:commutation_x_i}) between the $x_i$'s. Thus we have: 

\begin{proposition} \label{prop:RelationsSerres}
  The Serre relations
  \begin{gather}
    \sum_{r=0}^{1-a_{ij}}(-1)^r\qbinom{1-a_{ij}}{r}{q_i}\sym_z\left(x^\pm_i(z_1)\cdots
      x^\pm_i(z_r)x^\pm_j(w)x^\pm_i(z_{r+1})\cdots
      x^\pm_i(z_{1-a_{ij}})\right) \label{eq:serreannihilator}
  \end{gather}
  are in the annihilator ideal $\overline I_+$ of the weak Hopf
  pairing.  Moreover, the two-sided ideal generated by $I^1_+$ and
  (\ref{eq:serreannihilator}) is a weak Hopf ideal. A similar
  statement holds for the $x^-_i(z)$ ($i=1,\ldots,n-1$) and $\overline
  I_-$.
\end{proposition}

We will show that the Serre elements are quasi-primitive in
$\overline U_q(\b^+)/I^1_+$ and $\overline U_q(\b^-)/I^1_-$. We have
an induced weak Hopf pairing between these two algebras, and the
commutation relation (\ref{eq:commutation_x_i}) holds on top of
(\ref{eq:action_phi}).

Remark that a proof for $a_{ij} = -1, -2$ or $-3$ is sufficient.
Here we just give a straightforward computation in the case $a_{ij} =
-1$, covering in particular the case of $U_q(\hsl{n})$. The two
remaining cases were handled using a Computer Algebra System (Maple
V). Let us give some new notations. For $1\leq k \leq n$ and
$\epsilon\in\{0,1\}$ we put :
\[
X_{k,\epsilon} = 
\begin{cases}
  x^+_i(z_k)\ctensor 1 & \text{if $k<n$ and $\epsilon = 1$}, \\
  \phi_i(z_kq^{\frac{c_1}2})\ctensor x^+_i(z_kq^{c_1}) & \text{if $k<n$ and $\epsilon = 0$}, \\
  x^+_j(w)\ctensor 1 & \text{if $k=n$ and $\epsilon = 1$}, \\
  \phi_j(wq^{\frac{c_1}2})\ctensor x^+_j(wq^{c_1}) & \text{if $k=n$ and $\epsilon = 0$}. \\
\end{cases}
\] 
A similar definition can be given for $Y_{k,\epsilon}$ in $\overline
U_q(\b^-)/I^1_-\ctensor \overline U_q(\b^-)/I^1_-$. Now for
$\sigma\in\Sfrak_n$ and $\epsilon = (\epsilon_1\ldots\epsilon_n) \in
\{0,1\}^n$ we put $P^+_{\sigma,\epsilon} =
\prod_{k=1}^nX_{\sigma(k),\epsilon_{\sigma(k)}}$, $C^n_\sigma =
(-1)^{\sigma^{-1}(n)+1}\qbinom{n-1}{\sigma^{-1}(n)-1}{q_i}$ and
$S^+_{\epsilon} = \sum_{\sigma\in\Sfrak_n}
C^n_{\sigma}P^+_{\sigma,\epsilon}$. The coproduct of
(\ref{eq:serreannihilator}) in $\overline U_q(\b^+)/I^1_+\ctensor
\overline U_q(\b^+)/I^1_-$ is just
\[\sum_{\epsilon\in\{0,1\}^{2-a_{ij}}}S^+_{\epsilon}.\]

Let us give an example: take $\epsilon=(1,0,0)$. We get \[
P_{\text{id},\epsilon} = X_{1,1}X_{2,0}X_{3,0} = x_i^+(z_1)\tensor 1
\cdot \phi_i(z_2) \tensor x_i^+(z_2) \cdot \phi_j(w) \tensor
x_j^+(w). \] Thus, \[ P_{\text{id},\epsilon} = x_i^+(z_1) \phi_i(z_2) \phi_j(w)
\tensor x_i^+(z_2) x_j^+(w). \] 

We want to show that all the elements $S^+_\epsilon$ are zero but the
two extremal one ({\em i.e.} when all $\epsilon_i$ are $0$ or all
$\epsilon_i$ are $1$). Recall the remark we made on
page \pageref{remark-topology}. We observed that the topological
completion on the tensor product is just weak enough in order to allow
expressions of the form $g_{ij}(zw^{-1})x^+_i(z)\ctensor x^+_j(w)$ but
not $g_{ij}(zw^{-1})x^+_i(z)x^+_j(w)\ctensor 1$. To sketch the
forthcoming proof, let us just say that the completion allows us to
use commutation relation between $x^+_i$ and $x^+_j$ only when they
are not at the same side of the tensor product. The only cases where
such elements cannot be found are the two ``extremal'' one. This is
why the completion had to be carefully chosen. We need the following
lemmas:

\begin{lemma} \label{lemma:completion_implies_zero}
  Let $E=\sum_{n,m\in\Z}L_n\ctensor R_mz^{-n}w^{-m}$ be some
  generating series with $L_n$ and $R_m$ in $\overline
  U_q(\b^+)/I^1_+$.  Assume that for $n$ and $m$ fixed in $\Z$, the
  degree of $L_{n+k}\ctensor R_{m-k}$ goes to $+\infty$ (according to
  the tensor product filtration) when $k \rightarrow +\infty$. Then
  \hbox{$(az-bw)E = 0$} implies $E=0$ (where $a,b\in\C^*$).
\end{lemma}

Fix $n$ and $m$ in $\Z$. If $(az-bw)E = 0$ then we have
$aL_{n+1}\ctensor R_m = bL_n\ctensor R_{m+1}$ for all $n,m\in\Z$.
Going inductively, we get $(\frac{b}{a})^kL_{n+k}\ctensor R_{m-k} =
L_n\ctensor R_m$. As the left hand side degree goes to $+\infty$ when
$k\rightarrow +\infty$ and the right hand side is of constant finite
degree, we must have $L_n\ctensor R_m = 0$. As this is true for all
$n$, $m$ in $\Z$, we get $E = 0$. 

\begin{lemma} \label{lemma:serre_is_zero}
  In $\overline U_q(\b^+)/I^1_+ \ctensor \overline U_q(\b^+)/I^1_+$ we have:
  \begin{gather}
    (q^{-1}z_1-w)(q^{-1}z_2-w)S^+_{\{1,1,0\}} = 0 \,\notag\\
    (q^{2}z_1-z_2)(q^{-1}w-z_2)S^+_{\{1,0,1\}} = 0 \,\notag\\
    (q^{2}z_2-z_1)(q^{-1}w-z_1)S^+_{\{0,1,1\}} = 0 \,\notag\\
    (q^{2}z_1-z_2)(q^{-1}z_1-w)S^+_{\{1,0,0\}} = 0 \,\notag\\
    (q^{2}z_2-z_1)(q^{-1}z_2-w)S^+_{\{0,1,0\}} = 0 \,\notag\\
    (q^{-1}w-z_1)(q^{-1}w-z_2)S^+_{\{0,0,1\}} = 0 \,\notag
  \end{gather}
\end{lemma}

This lemma is the actual computation (which is difficult to handle by
hand for $a_{ij} = -2$ or $-3$). We will only show how to handle the
first relation. The remaining cases can be computed likewise. Using
only the action of $\phi_j(wq^{\frac{c}2})$ on $x^+_i(z_k)$ ({\em
  i.e.} the commutation relations
$(w-q^{-1}z_k)\phi_j(wq^{\frac{c}2})x^+_i(z_k) =
(q^{-1}w-z_k)x^+_i(z_k)\phi_j(wq^{\frac{c}2})$ ) we get :
\begin{align*}  
  (q^{-1}z_1-w) & (q^{-1}z_2-w)S^+_{\{1,1,0\}} = \\
  & -wq^{-2}(q-q^{-1}) \Big[
  (z_1-q^2z_2)x^+_i(z_1)x^+_i(z_2)\phi_j(wq^{\frac{c_1}2}) \\
  & \qquad - (z_2-q^2z_1)x^+_i(z_2)x^+_i(z_1)\phi_j(wq^{\frac{c_1}2})
  \Big] \ctensor x^+_j(cq^{c_1}).
\end{align*}  
Now using (\ref{eq:commutation_x_i}) (which is a valid relation
because we do this computation in $\overline U_q(\b^+)/I^1_+ \ctensor
\overline U_q(\b^+)/I^1_+$) we see the left operand of the tensor
product is zero.

Now remark that all the factor of the form $(az-bw)$ in lemma
\ref{lemma:serre_is_zero} are precisely those satisfying the condition
in lemma \ref{lemma:completion_implies_zero}. Thus we get that
$S^+_\epsilon = 0$ for $\epsilon \neq \{0,0,0\}$ and $\epsilon \neq
\{1,1,1\}$.  So the Serre relations are quasi-primitive in $\overline
U_q(\b^+)/I^1_+$.  Continuing as in proposition
\ref{prop:RelationsAnnihilator} we see the Serre relations are in the
annihilator ideal $\overline I_+$. A similar result holds for the
negative part.

Let $I^3_+$ and $I^3_-$ be the two sided ideals generated by all the
commutation relations we have seen so far. These are weak Hopf ideals,
and we denote by $U_q(\b^\pm)$ the quotient of $\overline U_q(\b^\pm)$
by $I^3_\pm$. The induced Hopf pairing leads to the weak quantum double
$\D(U_q(\b^+),U_q(\b^-))$. The new Hopf pairing will be noted as the
former, and let $I_+$ and $I_-$ be the annihilator ideals of this new
pairing. 


\subsection{An algebra morphism between $U_q(\hg)$ \\ and
  $\D(U_q(\b^+),U_q(\b^-))$} 

Remark that in $\D(U_q(\b^+), U_q(\b^-))$, the elements
$q^{\frac{c}2}\tensor q^{-\frac{c'}2}$ and $\phi_{i,0}\tensor
\psi_{i,0}$ are central and group like. Moreover, let $I$ be the
two-sided ideal generated by $q^{\frac{c}2} \tensor q^{-\frac{c'}2} -
1\tensor 1$ and $\phi_{i,0}\tensor \psi_{i,0} - 1\tensor 1$. We know
there is an algebra morphism from $U_q(\hg)$ onto $\D(U_q(\b^+),
U_q(\b^-)) / I$.  Actually, we have a little bit more :

\begin{proposition}
  There is an Hopf algebra isomorphism $\Phi$ between $U_q(\hg)$ and
  $\D(U_q(\b^+), U_q(\b^-)) / I$.
\end{proposition}

\proof Because of the remark above we know that $\Phi$ is onto. It
remains to show it is an isomorphism as a vector space. We have
$\D(U_q(\b^+), U_q(\b^-)) / I = U_+ \,\tensor\, \C[q^{\pm\frac{c}2}]
\,\tensor\, U_-$, where $U_+$ is generated by the $x^+_i(z)$ and
$\phi_i(z)$ for $i=1,\ldots,n-1$, and $U_-$ is constructed
accordingly. Now remark that any element of $U_q(\hg)$ can be assumed
to be in $U_+\cdot\C[q^{\pm\frac{c}2}]\cdot U_-$. Indeed, all the
commutation relations necessary for this operation involve finite sums
(actually we get infinite sums only when we commute $x^+_i(z)$ and
$x^+_j(w)$, but we do not have to do that here). Moreover, the above
decomposition is unique, which means $U_q(\hg) \simeq U_+ \,\tensor\,
\C[q^{\pm\frac{c}2}] \,\tensor\, U_-$. Rosso proved this assertion for
$U_q(\g)$ in \cite{Rosso:FiniteRepresentations}. His proof can be
easily extended to our case in a straightforward matter, so it will
not be done here. Finally we have $U_q(\hg) \simeq U_+ \tensor
\C[q^{\pm\frac{c}2}] \tensor U_- = \D(U_q(\b^+), U_q(\b^-)) / I$, and
the result follows.

\begin{proposition}
  The Hopf pairing between $U_q(\b^+)$ and $U_q(\b^-)$
  is non degenerate.
\end{proposition}

\proof Let $X^+$ be an element of $I_+$ of minimal degree $\alpha\in
Q_+$. In $\Delta(X^+)$, an element of degree $(\beta,\beta')$ with
$\beta$ and $\beta'$ non zero is in $I^+\ctensor U(\b^+) +
U(\b^+) \ctensor I^+$, so because the $\beta+\beta'$ is minimal this
element must be 0. Then $\Delta(X^+)$ is the sum of two elements of degree
$(\alpha,0)$ and $(0,\alpha)$. More specifically, we have $\Delta(X^+)
= X^+ \ctensor k + k' \ctensor X^+$ with $k$ and $k'$ in $U_q(\hh^+)$,
{\em i.e.} $X^+$ is quasi primitive.  Now, $X^+$ is 
quasi-commutative : if $Y\in U(\b^-)$ we have $X^+Y = \lambda
YX^+$, with $\lambda\in\C$. To see that, we compute $YX^+$ in the weak
quantum double, to get
\[YX^+ = \sum \hpairing{X^+_{(1)}, \S(Y_{(1)})}\hpairing{X^+_{(3)},
  Y_{(3)}} X^+_{(2)}Y_{(2)}.\]
But $X^+$ being in the kernel of the weak
Hopf pairing, only the terms involving the element of degree
$(0,\alpha,0)$ in $\Delta^2(X^+)$ can be non zero. If $Y$ is of degree
$\beta$ then because the weak Hopf pairing is $Q$-graded, we can eliminate
all the terms of the sum except those involving the element of degree
$(0,\beta,0)$ in $\Delta^2(Y)$. Combining this, we have \[YX^+ =
\hpairing{X^+_{(1)}, \S(Y_{(1)})}\hpairing{X^+_{(3)}, Y_{(3)}} X^+Y =
\lambda X^+Y.\]
Let $\Lambda$ in $Q_+$ and $L(\Lambda)$ be a highest weight module
with highest weight vector $v_\Lambda$. 
Actually we have $L(\Lambda) = U(\b^-)v_\Lambda$. Now
$X^+v_\Lambda = 0$, and so $X^+L(\Lambda) = 0$ because $X^+$ is
quasi-commutative. So $X^+$ is in the annihilator of every 
highest weight module. That is :
\[X^+ \in \bigcap_{\Lambda\in Q^+}\annu{U(\hg)}L(\Lambda).\]
But this intersection is zero according to
\cite{Joseph:QuantumGroupsPrimitiveIdeals}, which ends the proof.

We have some kind of ``functorial'' construction giving the whole
algebra from its positive part. Still, this Borel subalgebra
$U_q(\b^+)$ is given in term of generators and relations. Rosso
exhibited an interesting construction which we apply here, allowing us
to construct $U_q(\b^+)$ from a suitable Hopf bimodule.


\section{Quantum shuffle construction of $U_q(\b^+)$}

We begin by recalling some facts about quantum shuffle algebras,
following Marc Rosso's point of view. For more details, see
\cite{Rosso:QuantumShuffles}.


\subsection{Tensor algebra and cotensor coalgebra}
\label{sec:tensor-and-cotensor}

The following facts are due to Nichols. More details can be found in
\cite{Nichols:BialgebrasTypeOne}.  Let $H$ be a $k$-Hopf algebra over
a commutative field $k$ with invertible antipode $\S$, and $M$ a Hopf
bimodule over $H$ (\ie $M$ is a $H$-bimodule and a $H$-bicomodule).
Let $\delta_L$ and $\delta_R$ be the left and right coaction. 

We have two dual constructions over $H$ and $M$: the tensor algebra is
\[ T_H(M) = H \oplus \bigoplus_{n\geq 1}M^{\tensor_Hn},\] where
the multiplication is given by concatenation over $H$ for elements of
non zero degree, and left or right module action when one element is
in $H$. The tensor algebra as an universal property from which
$T_H(M)$ can be endowed with a Hopf algebra structure, where the
coproduct is the unique algebra map extending the coproduct on $H$ and
$\delta_L + \delta_R$ on $M$.

Dually, the cotensor coalgebra is defined as
\[ T_H^c(M) = H \oplus \bigoplus_{n\geq 1} M^{\cotensor_H n},\] where
$M\cotensor_H M$ is the kernel of $\delta_R\tensor\Id -
\Id\tensor\delta_L : M\tensor M \rightarrow M\tensor H\tensor M.$ The
coproduct is induced by the coproduct on $H$; on $M^{\cotensor_Hn}$
the component of bidegree $(i,j)$ of the coproduct is given by the
$\delta_L\tensor\Id$ when $i=0$, $\Id\tensor\delta_R$ when $j=0$, and
is induced by the map $(m_1\tensor\cdots\tensor m_n) \rightarrow
(m_1\tensor\cdots\tensor m_i) \tensor (m_{i+1}\tensor\cdots\tensor
m_n)$ otherwise. The counit is $\epsilon_H\circ\pi$ where $\pi$ is the
projection onto degree zero. Here again the cotensor coalgebra has an
universal property making it an Hopf algebra, where the multiplication
is the unique coalgebra map extending the usual multiplication on $H$
and defined by the module structure maps on degree $H\tensor M +
M\tensor H$.

Let $S_H(M)$ be the sub-Hopf algebra of $T_H^c(M)$ generated by $H$
and $M$. It is a Hopf bimodule, and it can be also obtained by the
following dual construction: the universal property on $T_H^c(M)$
allows us to define an unique Hopf algebra map $\Theta$ from $T_H(M)$ to
$T_H^c(M)$ induced by the natural isomorphisms on elements of degree
zero and one. Then $S_H(M)$ is the image of $\Theta$. If $I$ is the
kernel of $\Theta$, then we have also $S_H(M) \simeq T_H(M)/I$. 


\subsection{The quantum shuffle algebra}

In \cite{Rosso:QuantumShuffles}, Rosso brought Nichols work a step
further by considering a braiding introduced by Woronowicz in
\cite{Woronowicz:differential-calculus}. This braiding allows us to
give a precise description of the coproduct in $T_H(M)$ and, dually,
the product in $T_H^c(M)$. 

Let us consider the subspaces of left and right coinvariants of $M$:
\begin{align*}  
  M^L & = \{ m\in M \mid \delta_L(m) = 1\tensor m \} \\ 
  M^R & = \{ m\in M \mid \delta_R(m) = m\tensor 1. \} 
\end{align*} 
We know that $M$ is isomorphic to $M^R \tensor H$ with trivial right 
module and comodule structure.  Moreover, $M^R$ is a sub left comodule
of $M$, and a left module for the left adjoint action $h\cdot m = \sum
h_{(1)}m\S(h_{(2)})$. Similar properties hold for $M^L$ and the right
adjoint coaction.

The braiding $\sigma$ introduced by Woronowicz sends $M^R \tensor M^R$ to
himself.  It is defined by \[\sigma(m \tensor m') = \sum
m_{(-2)}m'\S(m_{(-1)}) \tensor m_{(0)},\] and it satisfies the
usual braid equation
\[ (\Id\tensor\sigma)(\sigma\tensor\Id)(\Id\tensor\sigma) =  
   (\sigma\tensor\Id)(\Id\tensor\sigma)(\sigma\tensor\Id). \] 

Let us note $V=M^R$. We denote by $\Sfrak_n$ the symmetric group
of $\{1, \ldots, n\}$, and by $s_i$ the transposition $(i,i+1)$
for $i = 1,\ldots, n-1$. If $p_1+\cdots+p_k = n$, we denote by
$\Sfrak_{p_1, \ldots, p_k}$ the set of $w\in\Sfrak_n$ such that
$w(1) < w(2) < \cdots < w(p_1)$, $w(p_1+1) < \cdots < w(p_1+p_2)$,
$\ldots$, and $w(p_1 + p_2 + \cdots + p_{k-1} + 1) < \cdots < w(p_1
+ \cdots + p_k)$. Such a $w$ is called a $(p_1,\ldots,p_k)$-shuffle.

The braid group $\Bfrak_n$ acts on $V^{\tensor n}$ in the usual
way: for $i = 1,\ldots, n-1$, we associate the $i$-th generator
$\sigma_i$ of $\Bfrak_n$ with $\Id_{V^{\tensor i-1}} \tensor
\sigma \tensor \Id_{V^{\tensor n-i-1}}$ on $V^{\tensor n}$. Let
$w$ be a permutation of the set $\{1, \ldots, n\}$. Then the lift
of $w$ in $\Bfrak_n$ is $T_w = \sigma_{i_1}\cdots\sigma_{i_k}$,
where $w = s_{i_1}\cdots s_{i_k}$ is a reduced expression of $w$.
The corresponding Hopf bimodule isomorphism in $V^{\tensor n}$ will
also be denoted by $T_w$.

\begin{proposition}
  Let $\sh$ be the product on $T(V)$ defined by:
  \[(x_1\tensor\cdots\tensor x_p)\sh(x_{p+1}\tensor\cdots \tensor x_n)
  = \sum_{w\in\Sfrak_{p,n}}T_w(x_1\tensor\cdots\tensor x_n),\] where
  $x_1,\ldots,x_n \in V$ and $\Sfrak_{p,n}$ is the set of
  $(p,n-p)$-shuffles. \\
  Then $\big( T(V), \sh \big )$ is an associative algebra.
\end{proposition}

The product $\sh$ is called a quantum shuffle product. This
construction is similar to the classical shuffle product, the usual
twist in $V\tensor V$ being replaced by $\sigma$. The algebra
$T(V)$ is then called a quantum shuffle algebra.

Remark that $T(V)$ is build on $V=M^R$. In order to have a more
general construction on $M=V \tensor H$, we consider $T(V) \tensor H$,
on which we put the following structure: it is an 
$H$-comodule, with $\delta_L$ given by the diagonal coaction of $H$ on
each $V^{\tensor n}$. We put the crossed product algebra structure on
$T(V) \tensor H$, with $H$ acting diagonally on $T(V)$. Finally, the
coalgebra structure is given by
\[\begin{split}
  \Delta(v_1\tensor\cdots\tensor v_n\tensor h) = \sum_{k=0}^n(v_1
  &\tensor \cdots \tensor v_k \tensor v_{k+1(-1)}\cdots v_{n(-1)}h_{(1)}
  ) \\
  & \tensor (v_{k+1(0)}\tensor\cdots\tensor v_{n(0)}\tensor h_{(2)}),
\end{split}\]
where $v_1,\ldots,v_n$ are in $V$, and $h\in H$. Those
structures are compatible, and $T(V) \tensor H$ becomes a Hopf
algebra.


\subsection{The quantum symmetric algebra}

Rosso showed that the cotensor coalgebra $T_H^c(M)$ is isomorphic to
$T(V) \otimes H$ as a right module and comodule. Furthermore, the
image of $T(V)$ in $T_H^c(M)$ is the subalgebra of right coinvariants
of $T_H^c(M)$.

Furthermore the sub Hopf algebra $S_H(M)$ of $T_H^c(M)$
generated by $H$ and $M$, is a Hopf bimodule, and it is
isomorphic to the crossed product of $H$ by $S_\sigma(V)$, where
$S_\sigma(V)$ is the subalgebra of $T(V)$ generated by $V$.
Actually, $S_\sigma(V)$ is isomorphic to the right coinvariants of
$S_H(M)$, via the previous isomorphism.

The Hopf algebra $S_H(M)$ is called a quantum symmetric algebra, and
is the main object of our study in the sequel of this paper.

 
\subsection{Construction of $S_H(M)$}

In order to apply Rosso's construction to $U_q(\b^+)$, we will
have to exhibit a sub Hopf algebra of it and put some Hopf bimodule
structure on $\tilde U_q(\b^+)$ over this subalgebra. What we gain by
using this functorial construction is a smaller number of generators
and relations (those necessary to describe the subalgebra), which
leads to simpler computations.

Let $H = U_q(\hh^+)$ be the sub Hopf algebra of $\complete{U_q(\b^+)}$
generated (algebraically, there is no need for a completion to make
$H$ an Hopf algebra) by $\phi_{i,l},\psi_{i,0}$ and
$q^{\pm\frac{c}2}$, for $i=1,\ldots,n-1,l \leq 0\}$. Let $V$ be the
subspace of $U_q(\hg)$ generated by $\X_{i,k}$, for $i=1,\ldots, n-1$
and $k \in \Z$. We would like to have a $H$-Hopf bimodule $V\tensor
H$, with a left action of $H$ reflecting relation
(\ref{eq:action_phi}) in definition \ref{def:Uq(g)}. But as this Hopf
bimodule will have to reflect the Hopf algebra structure of
$\complete{U_q(\b^+)}$ (see remarks below), we shall consider $V
\ctensor H$, the completion being similar to the one in
\ref{sec:gener-drinf-real}.

\begin{proposition}
  Let $M$ be the completed tensor product $V \ctensor H$. $M$ becomes
  an $H$-Hopf bimodule the following way: $V \ctensor H$ is a trivial
  right module and comodule. The left action of $H$ on $M$ is given by
  \begin{equation}
    \phi_{i,n}(x^+_{j,p} \ctensor \phi_{l,q}) = 
    \sum_{k\geq 0}g_k^{(ij)}x^+_{j,p-k} \ctensor q^{-k\frac{c}2}\phi_{i,n+k}\phi_{l,q}, \label{eq:actionHonM}
  \end{equation}
  and the left coaction of $H$ on $V$ by
  \[\delta_L(x^+_i(z)) = \phi_i(zq^{\frac{c_1}2}) \ctensor x^+_i(zq^{c_1}).\]
  The left coaction of $H$ on $M$ is then the diagonal coaction on $V
  \ctensor H$.
\end{proposition}

\remark{1} If we identify $x^+_i(z)\ctensor 1$ with $x^+_i(z)$ and $1
\ctensor \phi_i(z)$ with $\phi_(z)$ in $M$, then the left action of
$H$ on $M^R = V\ctensor 1$ can be written in the more satisfying form of
relation (\ref{eq:action_phi}).

\remark{2} If we compute $(\delta_L+\delta_R)(x^+_i(z))$, we recognize
the expression of $\Delta(x^+_i(z))$ in $\complete{U_q(\b^+)}$. This is
due to the fact that in the quantum symmetric algebra $S_H(M)$, the
coproduct is just $\delta_L+\delta_R$ on elements of degree $1$.

\proof We have to check two things. It is trivial to see that
$\delta_L$ and $\delta_R$ commutes, \ie $(\delta_L \tensor \Id)
\delta_R = (\Id \tensor \delta_R) \delta_L$. It remains to show that
(\ref{eq:actionHonM}) is a left action, and that $\delta_L$ and
$\delta_R$ are morphisms of $H$-bimodules. This can be done by a
direct calculation, but generating series are not well suited for
that. Therefore we have to use the following relations, which are
translations of all the former relations involving those series:
\begin{align}
  \phi_{i,n}\phi_{j,m} & = \phi_{j,m}\phi_{i,n} \,,\notag\\
  (x^+_{i,n}\ctensor \phi_{j,m})\phi_{k,p} & = x^+_{i,n}\ctensor \phi_{j,m}\phi_{k,p} \,,\notag\\
  \delta_R(x^+_{i,n}\ctensor 1) & = (x^+_{i,n} \ctensor 1) \ctensor 1 \,,\notag\\
  \delta_L(x^+_{i,n}) & = \sum_{k\leq 0}\phi_{i,k}q^{k\frac{c}2-nc} \ctensor x^+_{i,n-k} \,,\\
  \Delta(\phi_{i,n}) & = \sum_{r+s=n}\phi_{i,r}q^{-s\frac{c}2} \tensor \phi_{i,s}q^{r\frac{c}2} \,.
\end{align}
All verifications are now straightforward.

In the next part, we will compute the braiding associated to $M$.


\subsection{The braiding on $M^R \tensor M^R$} 

Recall that the braiding $\sigma$ on $M^R \tensor M^R$
is defined by \[\sigma(m \tensor m') = \sum m_{(-2)}m'\S(m_{(-1)}) \tensor
m_{(0)}.\] 

\begin{proposition}
  The braiding $\sigma$ is given by:
  \[\sigma(x^+_i(z)\ctensor x^+_j(w)) = g_{ij}(zw^{-1})
  x^+_j(w)\ctensor x^+_i(z).\]
\end{proposition}

\proof We want to compute $(\Id\ctensor\S\ctensor\Id)(\Delta\ctensor\Id)\delta_L(x^+_i(z))$. 
That is   
        
\begin{align*}         
  (\Id\ctensor\S\ctensor\Id) & (\Delta\ctensor\Id)\delta_L(x^+_i(z)) =       
  (\Id\ctensor\S\ctensor\Id)(\Delta\ctensor\Id)\left(\phi_i(zq^{\frac{c_1}2}) \otimes x^+_i(zq^{c_1})\right) \\        
  & = (\Id\ctensor\S\ctensor\Id)(\Delta\ctensor\Id)\Big(\sum_{n\leq 0,m\in\Z}\phi_{i,n}q^{-n\frac{c}2-mc}\ctensor x^+_{i,m}z^{-n}z^{-m}\Big) \\       
  & = (\Id\ctensor\S\ctensor\Id)\Big( \sum_{n\leq 0,m\in\Z \atop r+s=n}\phi_{i,r}q^{-s\frac{c}2 -n\frac{c}2 - mc}\ctensor\phi_{i,s}q^{r\frac{c}2 -n\frac{c}2 - mc} \ctensor x^+_{im}z^{-n}z^{-m} \Big) \\     
  & = \sum_{n\leq 0,m\in\Z \atop r+s=n} \phi_{i,r}q^{-s\frac{c}2 -n\frac{c}2 - mc}\ctensor\S(\phi_{i,s})q^{-r\frac{c}2 +n\frac{c}2 + mc}\ctensor x^+_{i,m}z^{-n}z^{-m} \\     
  & = \sum_{n\leq 0,m\in\Z \atop r+s=n} \phi_{i,r}q^{-sc -r\frac{c}2 - mc}\ctensor\S(\phi_{i,s})q^{s\frac{c}2 + mc}\ctensor x^+_{i,m}z^{-n}z^{-m}     
\end{align*}                   
Therefore the braiding becomes:   
\begin{align*}    
  \sigma(x^+_i(z)\ctensor x^+_j(w)) & = \sum_{n\leq 0,m\in\Z \atop r+s=n}\phi_{i,r}q^{-sc -r\frac{c}2 - mc} x^+_j(w)\S(\phi_{i,s})q^{s\frac{c}2 + mc} z^{-n}\ctensor x^+_{i,m}z^{-m}\\  
  & = \sum_{n\leq 0 \atop r+s=n}\phi_{i,r}q^{-(r+s)\frac{c}2} x^+_j(w)\S(\phi_{i,s}) z^{-n}\ctensor x^+_i(z) \\  
  & = \sum_{n\leq 0,m\in\Z}q^{-n\frac{c}2}\sum_{r+s=n}\phi_{i,r} x^+_{j,m}\S(\phi_{i,s})z^{-n}w^{-m}\ctensor x^+_i(z) \\  
  & = \sum_{n\leq 0,
    m\in\Z}g^{(ij)}_{-n} x^+_{j,n+m}z^{-n}w^{-m}\ctensor x^+_i(z)\qquad\text{({\em cf.} remark 1)} \\  
  & = g_{ij}(zw^{-1})x^+_j(w) \ctensor x^+_i(z).
\end{align*}                
 
Using the $H$-Hopf bimodule $M$ with the braiding described above, we
take a more precise look at the quantum symmetric algebra we get.
 

\subsection{An isomorphism between $S_H(M)$ and $U_q(\b^+)$}

Now we state the main result of this paper.

\begin{theorem}
  There is a Hopf algebra isomorphism between $U_q(\b^+)$ and $S_H(M)$.
\end{theorem}

\remark{3} Considering the definition of $H$ and $M$, the existence of
such an isomorphism is not a surprise. Almost all the work is already
done, and the interesting fact is the actual existence of a quantum
symmetric algebra which is isomorphic to $U_q(\b^+)$.

\proof There is an obvious map going from the associative algebra with 
unit $1$ and generators $\{x^+_{i,k},\, \phi_{i,l},\, \phi_{i,0}^{-1}
,\, q^{\pm\frac{c}2} | i = 1,\ldots,n-1, k\in\Z, l\leq 0\}$ to $S_H(M)$
with the shuffle product. We now have to check that the relations
defining $U_q(\b^+)$ are verified in $S_H(M)$ for the quantum
shuffle product. Relation (\ref{eq:phicommutation}) of
definition~\ref{def:Uq(g)} is true by construction of $H$ and we
already examined relation (\ref{eq:action_phi}) in remark 1. We get
the commutation relation (\ref{eq:commuteannihilator}) between
$\X_i(z)$ and $\X_j(w)$ as follows: 
\[\X_i(z)\sh\X_j(w) = \X_i(z)\ctensor \X_j(w) + g_{ij}(zw^{-1})\X_j(w)\ctensor\X_i(z),\]
while on the other side:
\[\X_j(w)\sh\X_i(z) = \X_j(w)\ctensor \X_i(z) +
g_{ji}(wz^{-1})\X_i(z)\ctensor\X_j(w).\]
Then we get the result by
using the relations $(z-q_{ij}w)g_{ij}(zw^{-1}) = (q_{ij}z-w)$ and 
$(q^{ji}z - w)g_{ji}(wz^{-1}) = (z - q_{ij}w)$. Now it remains to show
the quantum Serre relations (\ref{eq:serreannihilator}). This is done
using the same computation than in proposition
\ref{prop:RelationsSerres} when we wanted to show that the Serre
elements are quasi primitives. That is, we show that the Serre
relation multiplicated by suitable factors is zero. This is easily
done by applying the commutation relation
(\ref{eq:commuteannihilator}) between $\X_i(z)$ and $\X_j(w)$. Now,
because of the nature of the shuffle product any multiplicative
factor we used satisfy the condition in lemma
\ref{lemma:completion_implies_zero}. Thus the Serre relations are
zero. 

Now we have an algebra morphism $\Phi$ from $U_q(\b^+)$ to $S_H(M)$ .
It is easy to show it is actually a Hopf algebra morphism (see remark
2). This morphism is onto by construction. To achieve the proof of
the theorem it remains to show that the morphism is one to one. The
quantum symmetric algebra $S_H(M)$ is a graded algebra. The elements
of $H$ are given degree 0, and those of $V = M^R$ are given degree 1.
We then use the following lemma:

\begin{lemma}
  Let $x$ be an element of degree at least 2 in the kernel of the
  morphism $\Phi$ from $U_q(\b^+)$ to $S_H(M)$ . Then $x$ is 0.
\end{lemma}

This is done by induction. We know that elements of degree 0 and 1 are
not in the kernel of $\Phi$. Let us suppose that there are no element
of degree at most $n$ in the kernel, and let $x$ be an element of
degree $n+1$. Then we have
\[\Delta(x) = \delta_L(x) + \sum x_{(0)}\ctensor x_{(1)}+
\delta_R(x).\] But the kernel is a Hopf ideal and the elements
$x_{(0)}$ and $x_{(1)}$ are of degree at most $n$. So we get finally
\[\Delta(x) = \delta_L(x) + \delta_R(x),\] 
which means $\Delta(x)$ is in $H\ctensor M + M\ctensor H$. But we know
that elements of degree at least 2 with such a coproduct are in the
kernel of the weak Hopf pairing between $U_q(\b^+)$ and $U_q(\b^-)$.
Using the non degeneracy of this weak Hopf pairing, we get that those
elements are null in $U_q(\b^+)$.

\subsection{Work of B. Enriquez}

After the preparation of the preliminary version of this manuscript,
the work of Enriquez \cite{Enriquez:CorrelationFunctions} came to
our attention. 

Though the author discusses a shuffle algebra description of the
positive part of $U_q(\g)$, his approach is completely different from
the one being considered here. Enriquez describes some vanishing
conditions on the correlation functions of Drinfeld currents of the
positive nilpotent part of $U_q(\g)$. Those conditions are then used
to give an isomorphism between the positive part of $U_q(\g)$ and some
shuffle algebra construction, though his shuffle algebra seems to be
different from ours. In his framework the problem of topological
completion can be completely avoided because only highest weight
modules are considered.  Moreover, Enriquez does not consider the case
with non zero central part (\ie he assumes that $c=0$).



\end{document}